\documentclass[12pt]{article}

\usepackage{amsmath,amsfonts,amsthm,amssymb}
\usepackage{fullpage}
\usepackage{float}
\usepackage{comment}

\theoremstyle{plain}
\newtheorem{theorem}{Theorem}

\theoremstyle{definition}

\def\al{\alpha}

\theoremstyle{remark}

\title{Three Series for the Generalized Golden Mean}

\author{Kevin Hare \\
Pure Mathematics \\
University of Waterloo \\
Waterloo, ON  N2L 3G1 \\
Canada\\
{\tt kghare@uwaterloo.ca} \\
\ \\
Helmut Prodinger\\
Department of Mathematical Sciences\\
Stellenbosch University\\
7602 Stellenbosch\\
South Africa\\
{\tt hproding@sun.ac.za} \\
\ \\
Jeffrey Shallit \\
Computer Science  \\
University of Waterloo \\
Waterloo, ON  N2L 3G1 \\
Canada\\
{\tt shallit@cs.uwaterloo.ca} 
}

\begin{document}

\maketitle

\begin{abstract}
As is well-known, the ratio of adjacent Fibonacci numbers tends to
$\phi = (1+\sqrt{5}\,)/2$, and the ratio of adjacent Tribonacci numbers
(where each term is the sum of the three preceding numbers) tends to
the real root $\eta$ of $X^3 - X^2 - X - 1 = 0$.
Letting $\alpha_n$ denote
the corresponding ratio for the generalized
Fibonacci numbers, where each term is the sum of the $n$ preceding,
we obtain rapidly converging series for $\alpha_n$, $1/\alpha_n$, and
$1/(2-\alpha_n)$.
\end{abstract}

\section{Introduction}

The Fibonacci numbers are defined by the recurrence
$$ F_i = F_{i-1} + F_{i-2}$$
with initial values $F_0 = 0$ and $F_1 = 1$.
The well-known Binet formula (actually already known to de Moivre)
expresses $F_i$ as a 
linear combination of the zeroes $\phi \doteq 1.61803 > 0 > \hat{\phi} $
of the characteristic polynomial of the recurrence $X^2 - X - 1$:
$$ F_i = {{\phi^i - \hat{\phi}^i} \over {\phi - \hat{\phi}}} .$$
Here the number $\phi = {{\sqrt{5} + 1} \over 2}$ is popularly
referred to as the {\it golden mean} or {\it golden ratio}.

Similarly, the ``Tribonacci'' numbers (the name is apparently due
to Feinberg \cite{Feinberg}; also see \cite{Sharp}) are  defined by
$$ T_i = T_{i-1} + T_{i-2} + T_{i-3}$$
with initial values $T_0 = T_1 = 0$ and $T_2 = 1$.
Here we also have that $T_i$ is a linear combination
of $\eta_1^i, \eta_2^i, \eta_3^i$, where $\eta_1,\eta_2,\eta_3$
are the zeroes of the characteristic polynomial
$X^3 - X^2 - X - 1$; see, e.g., 
\cite{Spickerman82}.  Here 
$$\eta_1 = {1\over 3}\left(1 + \root3\of{19+3\sqrt{33}} +
\root3\of{19-3\sqrt{33}} \right)$$
is the only real zero and $\eta_1 \doteq 1.839$.

The ``Tetranacci'' (aka ``Tetrabonacci'', ``Quadranacci'') numbers
are defined analogously by
$$ A_i = A_{i-1} + A_{i-2} + A_{i-3} + A_{i-4}$$
with initial values $A_0 = A_1 = A_2 = 0$ and $A_3 = 1$.
Once again, the $A_i$ can be expressed as a linear
combination of the zeroes of the 
characteristic polynomial $X^4 - X^3 - X^2 - X - 1$;
see, for example \cite{Lin}.

More generally, we can define the generalized Fibonacci
sequence of order $n$ by
$$ G_i^{(n)} = G_{i-1}^{(n)} + \cdots + G_{i-n}^{(n)} $$
with appropriate initial terms.
Here the associated characteristic polynomial is
$X^n - X^{n-1} - \cdots - X - 1$.  As is well-known \cite{Miles,Miller}, this
polynomial has a single positive zero $\alpha_n$,
which is strictly between
$1$ and $2$.  (The other zeroes are discussed in 
\cite{Zhu&Grossman}.)  Table~\ref{table1} gives decimal approximations
of the first few dominant zeroes.
Furthermore, as Dresden has shown \cite[Theorem 2]{Dresden}, knowledge of
$\alpha_n$ suffices to compute the $i$'th generalized Fibonacci
number of order $n$.

\begin{table}[H]
\begin{center}
\begin{tabular}{r|c}
$n$ & $\alpha_n$ \\
\hline
2 & 1.61803398874989484820 \\
3 & 1.83928675521416113255 \\
4 & 1.92756197548292530426 \\
5 & 1.96594823664548533719 \\
6 & 1.98358284342432633039 \\
7 & 1.99196419660503502110 \\
8 & 1.99603117973541458982 \\
9 & 1.99802947026228669866 \\
10 &1.99901863271010113866  \\
\end{tabular}
\caption{Generalized golden means}
\end{center}
\label{table1}
\end{table}

It is natural to wonder how the generalized golden means
$\alpha_n$ behave as $n \rightarrow
\infty$.  Dubeau \cite{Dubeau} proved that
$(\alpha_n)_{n \geq 2}$ is an increasing sequence that converges
to $2$.
In fact, it is not hard to show, using the binomial theorem,
that 
$$ 2 - {1 \over {2^n - {n \over 2} - {{n^2} \over {2^n}}}} < \alpha_n <
2 - {1 \over {2^n - {n \over 2}}} $$
for $n \geq 2$; see \cite{Forsyth}.

In this paper,
we give three series that approximate $\alpha_n$, $1/\alpha_n$, and
$1/(2-\alpha_n)$ to any desired order.   Remarkably, all three have
similar forms.

\begin{theorem}
Let $n \geq 2$, and define $\alpha = \alpha_n$, the positive real zero of
$X^n - X^{n-1} - \cdots - X - 1$. Let $\beta = 1/\alpha$.    Then
\begin{itemize}
\item[(a)] 
    \begin{equation*}
    \beta=\frac12+\frac12\sum_{k\ge1}\frac1k\binom{k(n+1)}{k-1}\frac1{2^{k(n+1)}}.
    \end{equation*}
    \label{part1}
\item[(b)]
    \begin{equation*}
    \al=2-2\sum_{k\ge1}\frac1k\binom{k(n+1)-2}{k-1} \frac1{2^{k(n+1)}}.
    \end{equation*}
    \label{part2}
\item[(c)]
    \begin{equation*}
    {1 \over {2-\alpha}} = 2^n - {n \over 2} - 
     \frac12 \sum_{k\ge1}\frac1k \binom{k(n+1)}{k+1} \frac{1}{2^{k(n+1)}}
    \end{equation*}
    \label{part3} 
\end{itemize}
\label{main}
\end{theorem}
The proof is given in the next three sections.
Our main tool is the classical Lagrange inversion formula;
see, for example, \cite[\S A.6, p.\ 732]{Flajolet&Sedgewick:2009}:

\begin{theorem}
\label{thm2}
Let $\Phi(t)$ and $f(t)$ be formal power series with $\Phi(0) \not= 0$,
and suppose $t = z \Phi(t)$.
If $\Phi(0) \neq 0$, we can write $t = t(z)$ as a formal power series in $z$.
Then
\begin{itemize}
\item[(a)]
$[z^k] t   = {1 \over k} [t^{k-1}] (\Phi(t))^k$;

\item[(b)]
$[z^k] f(t) = {1 \over k} [t^{k-1}] f'(t) (\Phi(t))^k$;
\end{itemize}
where, as usual, $[z^k] t$ (resp., $[z^k]f(t)$) denotes the coefficient
of $z^k$ in the series for $t$ (resp., $f(t)$).
\end{theorem}

\section{A series for $\beta$}

In this section, we will prove Theorem~\ref{main} (a), namely:
\[
    \beta=\frac12+\frac12\sum_{k\ge1}\frac1k\binom{k(n+1)}{k-1}\frac1{2^{k(n+1)}}.
\]

\begin{proof}
From 
$$ \alpha^n = \alpha^{n-1} + \cdots + \alpha + 1 $$
we get
$$ (1-\alpha)\alpha^n =  1-\alpha^n$$
and hence
\begin{equation}
\al^{n+1}-2\al^n+1=0 .
\label{eq1}
\end{equation}
Recalling that $\beta = 1/\alpha$ we get
\begin{equation}
\beta=\frac12+\frac12\beta^{n+1}.
\label{eq2}
\end{equation}
Let $\Phi(t) = (t+\frac12)^{n+1}$ and 
\begin{equation}
t=z \Phi(t),
\label{eq3}
\end{equation}
as in the hypothesis of Theorem \ref{thm2}.
We notice that $t = \beta-\frac12$ and $z = \frac12$ is a solution to
Eq.~\eqref{eq3},
as shown in Eq.~\eqref{eq2}.
From the Lagrange inversion formula and the binomial theorem, we get
\begin{equation*}
[z^k]t=\frac1k[t^{k-1}]\Big(t+\frac12\Big)^{k(n+1)}=\frac1k\binom{k(n+1)}{k-1}\frac1{2^{k(n+1)+1-k}}.
\end{equation*}
So
\begin{equation*}
t=\sum_{k\ge1}\frac1k\binom{k(n+1)}{k-1}\frac1{2^{k(n+1)+1-k}}z^k .
\end{equation*}
In particular, at $z = \frac12$ and $t = \beta - \frac12$, we get
\begin{equation*}
\beta=\frac12+\frac12\sum_{k\ge1}\frac1k\binom{k(n+1)}{k-1}\frac1{2^{k(n+1)}} ,
\end{equation*}
as required.
\end{proof}

\section{A series for $\alpha$}

In this section, we will prove Theorem~\ref{main} (b), namely:
\[
    \al=2-2\sum_{k\ge1}\frac1k\binom{k(n+1)-2}{k-1} \frac1{2^{k(n+1)}}.
\]
This formula was previously discovered in 1998 by
Wolfram \cite[Theorem 3.9]{Wolfram}.

\begin{proof}
From \eqref{eq1} we get
$$ \al^{n}(\al-2)+1=0 $$
and so
\begin{equation}
2-\al=\al^{-n}.
\label{eq4}
\end{equation}
Let $\Phi(t) = (1-\frac t2)^{-n}$ and 
\begin{equation*}
t=z \Phi(t)
\end{equation*}
as in the hypothesis of Theorem~\ref{thm2}.
We observe that $t = 2-\al$ and $z = 2^{-n}$ is a solution, 
    as shown in Eq.~\eqref{eq4}.
Using the Lagrange inversion formula again, we find
\begin{equation*}
[z^k]t=\frac1k[t^{k-1}]\Big(1-\frac t2\Big)^{-kn}=\frac1k\binom{k(n+1)-2}{k-1}
\frac1{2^{k-1}}.
\end{equation*}
Therefore
\begin{equation*}
t=\sum_{k\ge1}\frac1k\binom{k(n+1)-2}{k-1}z^k 
\frac1{2^{k-1}}.
\end{equation*}
In particular, evaluating this at $t = 2-\al$ and $z = 2^{-n}$ gives
\begin{equation*}
2-\al=\sum_{k\ge1}\frac1k\binom{k(n+1)-2}{k-1}2^{-nk}
\frac1{2^{k-1}},
\end{equation*}
or
\begin{equation*}
\al=2-2\sum_{k\ge1}\frac1k\binom{k(n+1)-2}{k-1}
\frac1{2^{k(n+1)}}, 
\end{equation*}
giving us a series for $\alpha$.
\end{proof}

\section{A series for $1/(2-\alpha)$}

In this section we will prove Theorem~\ref{main} (c), namely:
    \begin{equation*}
    {1 \over {2-\alpha}} = 2^n - {n \over 2} - \frac12 \sum_{k\ge1}\frac1k \binom{k(n+1)}{k+1} \frac{1}{2^{k(n+1)}} .
    \end{equation*}

\begin{proof}
Define
\begin{align*}
S(z) =-\frac12\sum_{k\ge1}\frac1kz^k[t^{k+1}](1+t)^{k(n+1)}.
\end{align*}
At $z = 2^{-(n+1)}$, this gives
\begin{align*}
S\left(1/2^{n+1}\right)& = -\frac12\sum_{k\ge1}\frac1k\binom{k(n+1)}{k+1}\frac1{2^{k(n+1)}} .
\end{align*}
Hence it suffices to show that
\[S\left(1/2^{n+1}\right) = -2^n + {n \over 2} + \frac{1}{2-\alpha}. \]

We see from Eq.~\eqref{eq4} that
\begin{equation}
{2\over\alpha} -1 = \alpha^{-n-1} 
\label{eqn}
\end{equation}

Let $t = z \Phi(t)$ as before. 
Further let
\begin{equation*}
\Phi(t)=(1+t)^{n+1},\quad f'(t)=-\Phi^{-2}.
\end{equation*}
We see that $z = 1/2^{n+1}$ and $t = {2\over \alpha} -1$ is a solution
to $t = z \Phi(t)$ by Eq.~\eqref{eqn}.

To get a series for $1/(2-\alpha)$, we start from the Lagrange inversion formula, part (b), 
    to get
\begin{equation*}
f(t)=f(0)+\sum_{k\ge1}\frac1kz^k[t^{k-1}] (\Phi(t))^k f'(t).
\end{equation*}
Differentiating with respect to $z$ gives
\begin{equation*}
\frac{d}{dz}f(t)=\frac{dt}{dz}\cdot f'(t)=\sum_{k\ge1}z^{k-1}[t^{k-1}]
(\Phi(t))^k f'(t).
\end{equation*}
Using $t = z \Phi(t)$ we see that $\frac{dt}{dz} = \frac{\Phi(t)^2}{\Phi(t) - \Phi'(t)}$.
This gives us
\begin{align*}
\frac{\Phi^2}{\Phi-t\Phi'}\cdot f'(t)&=\sum_{k\ge1}z^{k-1}[t^{k-1}] (\Phi(t))^k f'(t)\\
&=[t^{0}]\Phi(t)f'(t)+
z^{1}[t^{1}] (\Phi(t))^2 f'(t)+\sum_{k\ge1}z^{k+1}[t^{k+1}] (\Phi(t)^k) 
(\Phi(t))^2 f'(t).
\end{align*}

Using the fact that $f'(t) = -\frac{1}{\Phi^2}$ we get
\begin{align*}
-\frac{1}{\Phi-t\Phi'}=-1-\sum_{k\ge1}z^{k+1}[t^{k+1}] (\Phi(t))^k.
\end{align*}
Observing that $S'(z) = -\frac12 \sum_{k\ge1}z^{k-1}[t^{k+1}](1+t)^{k(n+1)}$,
this simplifies to
\begin{equation*}
2z^2S'(z)=1-\frac{1}{\Phi-t\Phi'} .
\end{equation*}
Thus
\begin{equation*}
S'(z)=\frac1{2z^2}-\frac{1}{\Phi-t\Phi'}\frac{\Phi^2}{2t^2},
\end{equation*}
so
\begin{equation*}
S(z)=-\frac1{2z^1}-\int \frac{1}{\Phi-t\Phi'}\frac{\Phi^2}{2t^2} dz
=-\frac1{2z}-\int \frac{\Phi-t\Phi'}{\Phi^2}\frac{1}{\Phi-t\Phi'}\frac{\Phi^2}{2t^2} dt
\end{equation*}
and
\begin{equation*}
S(z)=-\frac1{2z}-\int dt\frac{1}{2t^2}=-\frac1{2z}+\frac{1}{2t}+C.
\end{equation*}
In order to compute the integration constant $C$, we note that $S(0)=0$. Then
\begin{equation*}
C=\frac12\lim_{z\to0}\Big[\frac1z-\frac1t\Big]=
\frac12\lim_{t\to0}\frac{\Phi-1}{t}=
\frac12\lim_{t\to0}\frac{(1+t)^{n+1}-1}{t}=\frac{n+1}{2}
\end{equation*}
and
\begin{equation*}
S(z)=-\frac1{2z}+\frac{1}{2t}+\frac{n+1}{2}.
\end{equation*}

\begin{comment}
Recalling that $z=2^{-(n+1)}$, set
\begin{equation*}
\al=\frac2{1+t} .
\end{equation*}
Then
\begin{equation*}
t=\frac2\al-1 = \Big(\frac{1+t}{2}\Big)^{n+1},
\end{equation*}
or
\begin{equation*}
\frac2\al-1=\al^{-n-1};
\end{equation*}
hence
\begin{equation*}
2-\al=\al^{-n},
\end{equation*}
which checks.
\end{comment}

Evaluating at $z = 1/2^{n+1}$ and $t = {2\over\alpha} - 1$ we have
\begin{equation*}
S(1/2^{n+1})=-\frac1{2z}+\frac{1}{2t}+\frac{n+1}{2}=-2^n+\frac{\al}{2(2-\al)}
+\frac{n+1}{2}=-2^n+\frac{1}{2-\al}
+\frac{n}{2},
\end{equation*}
as required. 
\end{proof}

\section{Speed of convergence}

The speed of convergence of the series in Theorem \ref{main}
    is determined by the individual terms in the sequence.
For example, consider the  series for $1/\alpha$:
    \begin{equation*}
    \beta=\frac12+\frac12\sum_{k\ge1}\frac1k\binom{k(n+1)}{k-1}\frac1{2^{k(n+1)}}.
    \end{equation*}
The convergence depends upon the speed of convergence of 
\[ f_1(k,n)/2^{k(n+1)} := \frac1k \binom{k(n+1)}{k-1} \frac{1}{2^{k(n+1)}}. \]
Similarly define
\[ f_2(k,n)/2^{k(n+1)} := \frac1k \binom{k(n+1)-2}{k-1} \frac{1}{2^{k(n+1)}}. \]
\[ f_3(k,n)/2^{k(n+1)} := \frac1k \binom{k(n+1)}{k+1} \frac{1}{2^{k(n+1)}} \]
based on the expansion of $\alpha$ and $1/(2-\alpha)$.

Notice that, by Stirling's approximation, we have
\begin{align*}
\lim_{k\to\infty} \log_2(f_1(k,n))/k & \doteq  
   \lim_{k\to\infty} \log_2(f_2(k,n))/k  \\ & \doteq 
   \lim_{k\to\infty} \log_2(f_3(k,n))/k \\ & \doteq  
   (n+1) \log (n+1) - n \log_2(n)  ,
\end{align*}
which, as $n \to \infty$, tends to 
    \[ \log_2(n+1) + \frac{1}{\log(2)} . \]

Thus, for example, when $n = 2$ (corresponding to the Fibonacci case),
we have 
$$\log_2 f_i (k, n) \sim (3 \log_2(3) - 2 \log_2(2)) k \sim (2.75489 \cdots)k.$$
Since each term of 
the summation is of the form $f_i(k,n)/2^{k(n+1)}$, in the case 
$n = 2$, the $k$'th term is approximately $2^{-.24511 k}$. 
Thus, for example,
1000 terms of the series are expected to give at least 73 correct digits;
in fact, it gives 77 or 78 depending on the series.
Here by digits of accuracy, we mean
  $\lfloor-\log_{10}| {\rm actual}-{\rm estimate}|\rfloor$, which is the number of 
  correct decimal digits after the decimal point.
See Table~\ref{tab:series} for a summation of various predictions versus actual accuracy.

\begin{table}[H]
\begin{center}
\begin{tabular}{llllll}
$n$ & $k$ & Predicted & Actual & Actual & Actual  \\
  &   & accuracy  & accuracy ($\alpha$) & accuracy ($1/\alpha$) & accuracy ($1/(2-\alpha)$) \\
\hline
2 & 100 & 7 & 10 & 10 & 9 \\
2 & 1000 & 73 & 78 & 78 & 77 \\
2 & 10000 & 737 & 744 & 743 & 743  \\ 
10 & 10 & 18 & 23 & 23 & 21 \\
10 & 100 & 185 & 192 & 191 & 190 \\
10 & 1000 &  1856 & 1864 & 1863 &  1862 \\
100 & 2 & 55 & 87 & 86 & 83 \\
100 & 10 & 279 & 311 & 311 & 307 \\
100 & 100 & 2796 & 2830 & 2829 & 2826 \\
\end{tabular}
\caption{Predicted and actual accuracy of truncated series}
\label{tab:series}
\end{center}
\end{table}

We notice that convergence is much much faster for larger $n$.

\bigskip

\noindent{\it Acknowledgment.}    The third author wishes to thank
J\"urgen Gerhard for his assistance with Maple and his suggestions about
the problem.

\end{document}